\documentclass[11pt]{article}

\usepackage{amssymb,amsmath,amsthm}

\theoremstyle{plain} 
\newtheorem*{thm}{Theorem}
\newtheorem*{pro}{Proposition}

\newtheorem*{cor}{Corollary}
\theoremstyle{definition} 
\newtheorem{*defn}{Definition}

\theoremstyle{remark}
\newtheorem{rem}{Remark}

\newcommand{\circk}{\mathop{K}\limits^\circ}
\newcommand{\n}{\noindent} 
\newcommand{\vp}{\varepsilon}

\newcommand{\cl}[1]{\mathcal{#1}} 
\newcommand{\bb}[1]{\mathbb{#1}}
\newcommand{\ovl}{\overline}

\begin{document}

\begin{center}
PLURISUBHARMONIC DEFINING FUNCTIONS, GOOD VECTOR FIELDS, AND EXACTNESS  
OF A CERTAIN ONE FORM\footnote{Research supported in part by NSF grant
DMS-9801539 and by the Erwin Schr\"odinger International Institute for
Mathematical Physics.}$^{,2}$
\end{center}
\bigskip

\begin{center} Emil J.~Straube and Marcel K.~Sucheston\footnote{This  paper
is based on joint work done prior to Marcel Sucheston's tragic death in April
2000.} \end{center} \vspace{.5in}

\n {\bf Abstract.} We show that the approaches to global regularity of the
$\bar\partial$-Neumann problem via the methods listed in the title are
equivalent when the conditions involved are suitably modified. These modified
conditions are also equivalent to one that is relevant in the context of Stein
neighborhood bases and Mergelyan type approximation.\vspace{.5in}

This paper is concerned with the relationship between some conditions, listed
in the title, that are known to imply global regularity of the
$\bar\partial$-Neumann problem. While these conditions are clearly related,
they are known not to be equivalent. We show that under certain natural
modifications they become equivalent. Interestingly, these (modified)
conditions are also equivalent to one that is relevant in a somewhat different
context: \ there should exist {\em conjugate\/} normal fields that are
(approximately) holomorphic in weakly pseudoconvex directions. Under favorable
circumstances, this leads (in addition to global regularity) to the existence
of Stein neighborhood bases and to Mergelyan type approximation theorems
([\ref{B-F}], [\ref{F-N}]; see Remark~\ref{rem1} below).

While these results are of interest from the general point of view of
understanding global regularity of the $\bar\partial$-Neumann problem,
concrete motivation came from our work in [\ref{S-S}], where we observed that
in the special situation considered there, the construction of the vector
fields having good commutation properties with $\bar\partial$ is equivalent
to the construction of a defining function plurisubharmonic at the infinite
type points of the boundary of the domain (see [\ref{S-S}], Remark~5; see also
Remark~2 below for further details). 

For background on the $\bar\partial$-Neumann problem, we refer the reader to
[\ref{F-K}], [\ref{B-S2}],  [\ref{C-S}].
Denote by $\rho$ a smooth $(C^\infty)$ defining function for $\Omega$. In
[\ref{B-S1}], Boas and the first author formulated a (necessarily) technical
condition in terms of a family of vector fields that have good commutation
properties with $\bar\partial$. There should exist a positive constant $C>0$
such that for every $\vp>0$, there exists a vector field $X_\vp$ of type (1,0)
whose coefficients are smooth in a neighborhood $U_\vp$ in $\bb C^n$ of the
set $K$ of boundary points of $\Omega$ of infinite type and such that 
\begin{equation}\label{eq1}
C^{-1} < X_\vp\rho < C\quad \text{on}\quad K,
\end{equation}
and
\begin{equation}\label{eq2}
|\partial\rho([X_\vp, \partial/\partial\bar z_j])| < \vp \quad \text{on}\quad
K,\qquad 1\le j \le n.
\end{equation}
We will say for short that  such a family of vector fields is {\em transverse
to $b\Omega$ ((\ref{eq1})) and commutes approximately with $\bar\partial$
((\ref{eq2})) at points of\/} $K$. It was shown in [\ref{B-S1}] that if
$b\Omega$ admits such a family, then the $\bar\partial$-Neumann operators
$N_q$ and the Bergman projections $P_{q}$, $0 \le q\le n$, are continuous in
Sobolev norms for $s\ge 0$. (Here, we consider the standard $\cl L^2$-Sobolev
spaces.) A detailed discussion of the ``vector field method'' may also be
found in [\ref{B-S2}]. We note that the condition we have formulated here is
slightly more stringent than what is required in [\ref{B-S1}], in that
$X_\vp\rho$ is  required to be real on $K$ (but not necessarily near $K$).
But this is satisfied in [\ref{B-S1}] and the other situations where the
vector field method works to establish regularity in $W^s$ for all $s\ge 0$,
see [\ref{B-S3}], [\ref{B-S1}] and the recent [\ref{S}] and [\ref{S-S}]. It is
not clear at present how much more stringent this condition actually is.

It was already noted in [\ref{B-S3}] that in (\ref{eq2}), it suffices to
consider commutators with vector fields in Levi null directions: \ the
commutators in the remaining directions can be adjusted by modifying $X_\vp$
in complex tangential directions (see [\ref{B-S3}], proof of the lemma for
details). On the other hand, complex tangential components do not contribute
to the normal (1,0)-component of commutators with fields in Levi null
directions (this is a consequence of pseudoconvexity). Accordingly, we are led
to consider fields which are (real) multiples of the normal to the boundary. 

For such a field, it is immaterial whether we consider the normal
(1,0)-component of commutators with type (1,0) or type (0,1) derivatives
(in complex tangential  directions), equivalently, whether we consider the
normal or its conjugate. Doing the latter has the advantage that complex
tangential components of type (0,1) derivatives in Levi null directions are
automatically zero (details will be provided below). We are thus led to the
following notion. We say that $b\Omega$ admits  {\em a family of conjugate
normals which are approximately holomorphic in weakly pseudoconvex
directions\/} if there is a constant $C>0$ such that for all $\vp>0$, there
exists a vector field $N_\vp = e^{h_\vp} \sum\limits^N_{j=1}
\frac{\partial\rho}{\partial\bar z_j} \frac\partial{\partial z_j}$ such that
\begin{equation}\label{eq3} -C < h_\vp < C, \end{equation}
and
\begin{equation}\label{eq4}
\ovl Y(\ovl N_\vp)(z) = O(\vp),\quad Y(z)\in N(z),\qquad z\in K,
\end{equation}
where $N(z)$ denotes the null space of the Levi form at $z, Y$ is a (local)
section of $T^{1,0}(b\Omega)$ of unit length, and $\ovl Y$ {\em acts
componentwise on\/}  $\ovl N_\vp$.

When computing (normal components of) commutators, as in (\ref{eq2}), there is
a 1-form that arises naturally. Denote by $\eta$ a purely imaginary,
non-vanishing smooth one-form on $b\Omega$ that is zero on the complex tangent
space and its conjugate. Let $T$ be the purely imaginary tangential vector
field on $b\Omega$ that is orthogonal to the complex tangent space and its
conjugate (in the metric induced by $\bb C^n$) and that satisfies $\eta(T)
\equiv 1$ on $b\Omega$. Set $\alpha := -\cl L_T\eta$, that is, $\alpha$ is
minus the Lie derivative of $\eta$ in the direction of $T$. The form $\alpha$
was introduced into the literature by D'Angelo [\ref{DA1}], [\ref{DA2}]. The
relevance of $\alpha$ in the context of global regularity was recognized in
[\ref{B-S1}]. A crucial property of $\alpha$  is the following
closedness property (this property hinges on pseudoconvexity): \ the
differential $d\alpha$ restricted (pointwise) to the null space of the Levi
form vanishes. We refer the reader to [\ref{B-S1}], section~2, for this and
other properties of $\alpha$ (see also [\ref{B-S2}], pp.~97--98, [\ref{S-S}]).

We say that $\alpha$ {\em is approximately exact on the null space of the Levi
form\/} if there exists a constant $C>0$ such that for all $\vp>0$ there exists
a smooth real-valued function $h_\vp$ in a neighborhood $U_\vp$ (in $b\Omega$)
of $K$ such that on $K$
\begin{equation}\label{eq5}
1/C \le h_\vp \le C,
\end{equation}
and
\begin{equation}\label{eq6}
dh_\vp/N(p) = \alpha/N(p) + O(\vp).
\end{equation}
Here, $O(\vp)$ denotes a 1-form that satisfies $|O(\vp)(X)| \le \text{const. }
\vp|X|$. Note that although $\alpha$ depends on the choice of $\eta$,
whether or not $\alpha$ is approximately exact on the null space of the Levi
form does not:\ direct computation shows that if $\tilde\eta = e^g\eta$
for some smooth function $g$, then the corresponding form $\tilde\alpha$
differs from $\alpha$, on complex tangent vectors, by the differential of $g$.

The remaining condition alluded to in the title is as follows. We say that
$\Omega$ admits a defining function that is plurisubharmonic at the boundary
([\ref{B-S3}]), if there exists some smooth defining function $\rho$ whose
complex Hessian is positive semi-definite at all boundary points. This
condition is slightly more stringent than pseudoconvexity, which requires
positive definiteness only on the complex tangent space (rather than in all
directions). It was shown in [\ref{B-S3}] that if $\Omega$ admits a defining
function that is plurisubharmonic at the boundary (near the points of infinite
type is sufficient), then $b\Omega$ admits a family of vector fields
transversal to the boundary and commuting approximately with $\bar\partial$
(and the Bergman projection and the $\bar\partial$-Neumann operator are
regular in Sobolev norms). It was pointed out in [\ref{B-S1}], Remark~3, that
the existence of such a family of vector fields is actually a weaker property
than the existence of a defining function that is plurisubharmonic at the
boundary: \ examples may be obtained by considering domains that have as
suitable lower dimensional sections domains not admitting (even) a (local)
plurisubharmonic defining function. (See [\ref{B-S1}], Remark~3 for
references.) However, the proof in [\ref{B-S3}] does not use the full force of
the positive definiteness of the complex Hessian $L_\rho$ at boundary points;
it suffices to have this on the span of the null space of the Levi form and the
complex normal to the boundary.  On this span, however, it amounts to the
same to assume that the complex Hessian is zero. Indeed, if $\tilde\rho =
g\rho$, then the complex Hessian $L_{\tilde\rho}$ of $\tilde\rho$ satisfies
(on the boundary) $L_{\tilde\rho}(X,\ovl Y)(P) =$ $g(P) L_\rho(X,\ovl Y)(P) +
X\rho(P) \ovl Yg(P) + \ovl Y\rho(P) Xg(P)$.
Now if $Y(P) \in N(P)$, then the semidefiniteness of $L_\rho$ on the span of
$N(P)$ and the complex normal at $P$ implies that $L_\rho(X,\ovl Y)(P) = 0$
for $X$ in this span (it is only for this conclusion that plurisubharmonicity
of $\rho$ was used in [\ref{B-S3}]). Consequently, if we choose $g \equiv 1$ on
$b\Omega$, then also $L_{\tilde\rho}(X,\ovl Y)(P) = 0$. Extending $g$ from the
boundary so that the (real) normal derivative equals minus $L_\rho(Z,\overline
Z)$,
where $Z$ is the complex normal (normalized so that $Z\rho\equiv 1$ on
$b\Omega$) gives that in addition also $L_{\tilde\rho}(Z,\overline Z) = 0$.
Thus $L_{\tilde\rho}(P)$ is zero on the span of $N(P)$ and the  complex normal.
Finally, as with the other conditions, we only need this condition to be
satisfied approximately.
Accordingly, we say that $\Omega$ admits {\em a family of essentially
pluriharmonic defining functions\/} if there exists a constant $C>0$ such that
for all $\vp>0$ there exists a $(C^\infty)$ defining function $\rho_\vp$ for
$\Omega$ satisfying
\begin{equation}\label{eq7}
1/C\le |\nabla\rho_\vp| \le C,
\end{equation}
and
\begin{equation}\label{eq8}
\left|\sum_{j,k} \frac{\partial^2\rho_\vp(P)}{\partial z_j\partial \bar z_k}
w_j\ovl w_k\right| \le O(\vp)|w|^2 \qquad \forall~w\in \text{span}_{\bb C}\{N(P), L_n(P)\}
\end{equation}
for all boundary points $P$ in $K$. Here ${\rm span}_{\bb C}$
denotes the linear span over $\bb C$, and $L_n$ is the complex normal (see
below for the normalization we will use).

\begin{thm}
Let $\Omega$ be a smooth bounded pseudoconvex domain in $\bb C^n$. The
following are equivalent:
\begin{itemize}
\item[(i)] $\Omega$ admits a family of essentially pluriharmonic defining
functions
\item[(ii)] $\Omega$ admits a family of conjugate normals which are
approximately holomorphic in weakly pseudoconvex directions
\item[(iii)] $\Omega$ admits a transversal family of vector fields which
commute approximately with $\bar\partial$
\item[(iv)] the form $\alpha$ associated to some choice of $\eta$ (hence to
any choice) is approximately exact on the null space of the Levi form.
\end{itemize}
\end{thm}

\bigskip

\begin{rem}
We emphasize again that condition (i) in the theorem is indeed a
generalization of the notion of pluri{\em sub}harmonic defining function, as
explained in the discussion preceding the definition of a family of
essentially pluriharmonic defining functions. 
\end{rem}

\bigskip

We begin with the equivalence of (iii) and (iv). That (iv) implies (iii) was
observed by Boas and the first author in [\ref{B-S2}] (see the discussion at
the end of section~6, the ideas are from [\ref{B-S1}], [\ref{B-S3}]); we will
recall the main points here. Also, a version of (iii)~$\Rightarrow$~(iv), in
the case when the infinite type points are contained in submanifolds of the
boundary of a certain kind, was pointed out in [\ref{B-S1}], Remark~5.

It will be convenient to choose $\eta := \partial\rho-\bar\partial \rho$,
where $\rho$ is a smooth defining function for $\Omega$, and $T := L_n-\bar
L_n$, where $L_n := \frac2{|\nabla\rho|^2} \sum\limits^n_{i=1}
\frac{\partial\rho}{\partial\bar z_j} \frac\partial{\partial z_j}$. If $Y$ is
a local section of $T^{1,0}(b\Omega)$, then the definition of the Lie
derivative, the fact that $\eta(\ovl Y) \equiv 0$, and the fact that
$(\partial \rho + \bar\partial\rho)([L_n-\bar L_n, \ovl Y]) = 0$ give
\begin{equation}\label{eq9}
\alpha(\ovl Y) = 2\partial\rho([L_n,\ovl Y])
\end{equation}
(compare [\ref{DA2}], p.~92, [\ref{B-S1}], p.~231). If $L_1,\ldots, L_{n-1}$
denote local sections of $T^{1,0}(b\Omega)$ that span $T^{1,0}(b\Omega)$
(locally), then the fields $X_\vp$ in (iii) can be written (locally) as
\begin{equation}\label{eq10}
X_\vp = e^{h_\vp}L_n + \sum^{n-1}_{j=1} a^\vp_jL_j,
\end{equation}
with smooth functions $a^\vp_j$ $(1\le j \le n-1)$ and $h_\vp$. Computing
commutators with $\bar L_k$, $1\le k \le n-1$, and taking normal
(1,0)-components gives (keep in mind that with the normalization above,
$(\partial\rho-\bar\partial\rho)(L_n-\bar L_n) = 1$, but $\partial\rho(L_n) =
\frac12$)
\begin{align}
\partial\rho([X_\vp,\bar L_k]) &= \left(-\frac12 \bar L_k h_\vp + \partial\rho
([L_n, \bar L_k])\right) e^{h_\vp} + \sum^{n-1}_{j=1} a^\vp_j
\partial\rho([L_j,\bar L_k])\nonumber\\
\label{eq11}
&= (-dh_\vp (\bar L_k) + \alpha(\bar L_k)) \frac{e^{h_\vp}}2 +
\sum^{n-1}_{j=1} a^\vp_j \partial\rho ([L_j,\bar L_k]).
\end{align}
Now let $P \in K$ (i.e.\ $P$ is a point of infinite type),
and let $L_k(P) \in N(P)$. Then $\partial\rho([L_j,\bar L_k])(P)\break = 0$, by
pseudoconvexity of $b\Omega$. (A mixed term in a positive semidefinite
Hermitian form vanishes if one of the entries is a null direction of the
quadratic form.) Consequently
\begin{equation}\label{eq12}
\partial\rho([X_\vp,\bar L_k])(P) = (-dh_\vp(\bar L_k)(P) + \alpha(\bar
L_k)(P)) \frac{e^{h_\vp(P)}}2.
\end{equation}
Taking into account that both $h_\vp$ and $\alpha$ are real, (\ref{eq12})
shows that (iii) implies (iv). For the converse implication, fix $P\in K$. We
may assume that $L_1,\ldots, L_{n-1}$ are orthonormal and that they
diagonalize the Levi form {\em at\/} $P$, and that $L_1(P),\ldots, L_m(P)$ span
$N(P)$ (for some $m$ with $1\le m \le n-1$). For $m+1 \le j \le n-1$, set
\begin{equation}\label{eq13}
a^\vp_j := \left(\frac12 \bar L_j h_\vp(P) - \partial\rho ([L_n,\bar L_j])(P)
\right)e^{h_\vp(P)}/\partial\rho([L_j,\bar L_j])(P).
\end{equation}
The field $X^p_\vp$, defined by $X^p_\vp  := e^{h_\vp} L_n +
\sum\limits^{n-1}_{j=m+1} a^\vp_j L_j$, (where $h_\vp$ comes from (iv), and
satisfies (\ref{eq5}), (\ref{eq6})) then satisfies
\begin{equation}\label{eq14}
|\partial\rho([X^p_\vp, \bar L_j])(P)| \le \widetilde C\vp,\qquad 1\le j \le
n-1. \end{equation}
From here on, the argument is exactly the same as that in the proof of the
lemma in [\ref{B-S3}], pp.~85--86: (\ref{eq14}) extends by continuity into a
neighborhood of $P$, and patching finitely many of these locally defined
fields via a partition of unity ($K$ is compact) yields a field $X_\vp$ defined
in a neighborhood of $K$ in $b\Omega$ that has the required commutation
properties (i.e.\ (\ref{eq1}) and (\ref{eq2})) with sections of
$T^{(0,1)}(b\Omega)$ (note that the function $h_\vp$ is defined globally, so
that the normal component $e^{h_\vp}L_n$ is defined globally; consequently,
terms in the commutators arising from derivatives of the cut-off functions are
complex tangential and vanish when $\partial\rho$ is applied); the field
obtained in this way can be corrected by a field identically zero on $b\Omega$
to accommodate commutators with $\bar L_n$.

To see that (i) implies (iii), let $\rho_\vp = e^{h_\vp}\rho$ (thus defining
$h_\vp$), $P\in K$, and $\bar L_k(P)\in T^{0,1}_{b\Omega}(P)$. Then
\begin{align}
e^{h_\vp(P)} \partial\rho([e^{-h_\vp}L_n,\bar L_k])(P) &= \sum_j \bar L_k
\left(e^{-h_\vp} \frac{2}{|\nabla\rho|^2} \frac{\partial\rho}{\partial\bar
z_j}\right)(P) \frac{\partial\rho}{\partial z_j} (P) e^{h_\vp(P)}\nonumber\\
&= -\sum_j e^{-h_\vp(P)} \frac{2}{|\nabla\rho(P)|^2}
\frac{\partial\rho}{\partial \bar z_j} (P) \bar L_k\left(\frac{\partial
\rho}{\partial z_j} e^{h_\vp}\right)(P),\nonumber\\
\label{eq15}
&= -e^{-h_\vp} L_{\rho_\vp}(L_n, \bar L_k)(P),
\end{align}
where $L\rho_\vp$ denotes the complex Hessian of $\rho_\vp$.
We have used in (\ref{eq15}) that 
$\sum\limits_j \frac2{|\nabla\rho|^2} \frac{\partial\rho}{\partial\bar z_j}
\frac{\partial\rho}{\partial z_j} \equiv 1/2$, that
$\frac{\partial\rho}{\partial z_j} e^{h_\vp} = \frac\partial{\partial z_j}
(\rho e^{h_\vp})$ on $b\Omega$, and that $\bar L_k$ is tangential. If now
$L_k(P)\in N(P)$, we obtain from (\ref{eq8}) by polarization (in view of the
uniform bounds on $|\nabla\rho_\vp| \approx e^{h_\vp}$ given by (\ref{eq7}))
\begin{equation}\label{eq16}
|\partial\rho([e^{-h_\vp}L_n, \bar L_k])(P)| = O(\vp).
\end{equation}
From here on, the construction of the family of vector fields required in
(iii) proceeds as in the proof above of the implication
(iv)~$\Rightarrow$~(iii) (which is, as we pointed out, as in [\ref{B-S3}]).

Conversely, if (iii) is satisfied, then, setting $X_\vp = e^{h_\vp} L_n +$
complex tangential terms, we have (as in (\ref{eq11}), (\ref{eq12})),
\begin{equation}\label{eq17}
\partial\rho([X_\vp, \bar L_k])(P) = \partial\rho([e^{h_\vp} L_n,\bar L_k]) (P)
\end{equation}
for $L_k(P) \in N(P)$. (\ref{eq15}) gives (replacing $h_\vp$ by $-h_\vp$) for
$\rho_\vp := e^{-h_\vp}\rho$
\begin{equation}\label{eq18}
|L_{\rho_\vp}(L_n,\bar L_k)(P)| = O(\vp).
\end{equation}
By the discussion immediately preceding the definition of a family of
essentially pluriharmonic defining functions, (\ref{eq18}) suffices to obtain
such a family.

(\ref{eq15}) also shows that (ii) $\Rightarrow$ (iii) (again using that
commutators in directions not in the null space of the Levi form can be
adjusted by adding suitable complex tangential terms, as in the proof that
(iv)~$\Rightarrow$~(iii)).

Going in the other direction, we have with $X_\vp = e^{-h_\vp} L_n +$ complex
tangential terms, and $L_k(P)\in N(P)$ \begin{align}
\partial\rho([X_\vp, \bar L_k])(P) &= \partial\rho([e^{-h_\vp} L_n, \bar
L_k])(P)\nonumber\\
\label{eq19}
&= O(\vp).
\end{align}
Combining (\ref{eq19}) with (\ref{eq15}) gives, if we set $N_\vp := e^{h_\vp}
\sum\limits^n_{j=1} \frac{\partial\rho}{\partial \bar z_j}
\frac\partial{\partial z_j}$, that $\bar L_k(\ovl N_\vp)(P)$ (where $\bar L_k$
acts componentwise as in (\ref{eq4})) has inner product with $\bar L_n(P)$
that is $O(\vp)$. The inner products with $\bar L_1(P),\ldots,
\bar L_{n-1}(P)$ are zero solely by virtue of the fact that $L_k(P)\in N(P)$,
regardless of $h_\vp$. Indeed, fix $j\in \{1,\ldots, n-1\}$ and let $L_j(P)
= (\zeta_1,\ldots, \zeta_n)$, $L_k(P) = (w_1,\ldots, w_n)$. Then
\begin{align}
\langle \bar L_k(\ovl N_\vp)(P), \bar L_j(P)\rangle &=\sum_{j,k} \zeta_j \ovl
w_k \frac\partial{\partial\bar z_k} \left(e^{h_\vp}
\frac{\partial\rho}{\partial z_j}\right)(P)\nonumber\\
\label{eq20}
&= e^{h_\vp} \sum_{j,k} \frac{\partial^2\rho(P)}{\partial z_j\partial \bar
z_k} \zeta_j \ovl w_k\\
&= e^{h_\vp} L_\rho(L_j,\bar L_k)(P) = 0.\nonumber
\end{align}
We have used in the second equality that $\sum\limits_j
\frac{\partial\rho}{\partial z_j}(P) \zeta_j = 0$, and that $\Omega$ is
pseudoconvex and $L_k(P)\in N(P)$ in the last equality. Because
$\{\bar L_1(P), \bar L_2(P),\ldots, \bar L_{n-1}(P), \bar L_n(P)\}$ is a
basis (over $\bb C$) of $\bb C^n$, we obtain that $\bar L_k(\ovl N_\vp)(P)$ is
$O(\vp)$.

This concludes the proof of the theorem.

\bigskip

\begin{rem}\label{rem1}
As noted in the introduction, the existence of a family of conjugate normals
with suitable holomorphicity properties can lead to the existence of Stein
neighborhood bases and Mergelyan type approximation theorems. This follows from
[\ref{B-F}]. We briefly discuss one such instance here; it arises from our
work in [\ref{S-S}].

We consider smoothly bounded pseudoconvex domains in $\bb C^2$ whose set $K$
of infinite type boundary points  has a smooth boundary, as
 a subset of $b\Omega$. The interior of $K$ is then foliated by
1-dimensional complex manifolds. This foliation is usually referred to as the
Levi foliation of $\circk$. For generic such $K$, the conditions in the
theorem above are satisfied; in fact the conjugate normals can be taken to
be holomorphic along the leaves of the Levi foliation of $K$, that is, they
are $CR$-functions on $K$. 

$\Gamma$, the boundary of $K$, is a 2-dimensional
surface sitting inside $b\Omega$, and complex tangents occur precisely at
points of $\Gamma$ where the tangent space of $\Gamma$ coincides with the
complex tangent space to $b\Omega$. Recall that a generic complex tangency is
one that is either elliptic or hyperbolic. (This goes back to Bishop's paper
[\ref{B}]; see the introduction of [\ref{F}] for a thorough discussion of these
matters).

For terminology from foliation theory, in particular for the notion of
(infinitesimal) holonomy, the reader may consult [\ref{C-C}].

 Combining our
ideas from [\ref{S-S}] with ideas from [\ref{B-F}] yields the following
result, which may be viewed, in the case of $\bb C^2$, as a strengthening of
some aspects of Theorem~6.3 of [\ref{B-F}] (their theorem addresses the
situation in $\bb C^n$ for general $n$, however).

\begin{pro}
a)~~Let $\Omega$ be a smooth bounded pseudoconvex domain in $\bb C^2$. Suppose
that the set $K$ of weakly pseudoconvex boundary points is smoothly bounded
(in $b\Omega$) and that its boundary $\Gamma$ has only isolated generic complex
tangencies. Assume that the two leaves that meet at a hyperbolic point of
$\Gamma$ have no other hyperbolic points in their closure (in $K$). If each
leaf of the Levi foliation is closed (in $\circk$) and has trivial
infinitesimal holonomy, then there exists a conjugate normal field which is
$CR$ on $K$.

\n b)~~If in addition $K$ is uniformly $H$-convex, then there exists a
holomorphic vector field in a neighborhood of $K$ (in $\bb C^n$) that is
transverse to $b\Omega$ near $K$.
\end{pro}

It is part of the assumption in a) that the two local leaves that meet at a
hyperbolic point are globally distinct; see the discussion in [\ref{S-S}].
Note that if the leaves of the Levi foliation are assumed simply connected
(i.e.\ they are the analytic discs), then the (infinitesimal) holonomy is
automatically trivial. Recall that $K$ is uniformly $H$-convex if it admits a
Stein neighborhood basis $\{U_j\}^\infty_{j=1}$ of open pseudoconvex sets such
that for some constant $c>0$, $\{z\mid \text{dist}(z,K) < \frac1{cj}\}
\subseteq U_j$ $\subseteq\{z\mid \text{dist}(z,K)<\frac{c}j\}$.

\begin{cor}
Under the assumptions of the proposition, part b), we have
\begin{itemize}
\item[a)] $\ovl\Omega$ admits a Stein neighborhood basis
\item[b)] Functions analytic in $\Omega$ and continuous on $\ovl\Omega$ can be
approximated uniformly on $\ovl\Omega$ by functions holomorphic in some
neighborhood of $\ovl\Omega$.
\end{itemize}
\end{cor}

The corollary follows directly from the proposition, part b), and [\ref{B-F}],
\S 7, in particular Lemma~7.3 (for a)); and [\ref{F-N}], Theorem~1 (for b)).

To prove the proposition, we note that part a) comes from [\ref{S-S}],
Theorem~1, and the theorem above: \ Theorem~1 in [\ref{S-S}] gives a family of
vector fields in a neighborhood of $K$ that commute approximately with
$\bar\partial$, and our theorem above then gives a family of conjugate normals
which are approximately holomorphic in weakly pseudoconvex directions (i.e.\
along the leaves of the Levi foliation of $\circk$). Inspection of the proofs
(both in [\ref{S-S}] and in the theorem above) shows that actually these
conjugate normal fields can be taken to be exactly holomorphic on the leaves,
i.e.\ $CR$ on $K$.

The proof of b) now follows entirely by the arguments in [\ref{B-F}]. If $K$
is uniformly $H$-convex, $CR$-functions in $K$ can be approximated uniformly
on $K$ by functions holomorphic in a neighborhood of $K$ (see e.g.\
[\ref{B-F}], proof of Proposition~6.2). Denoting the algebra of these
functions by $A(K)$, we may furthermore invoke Theorem~2.12 in [\ref{R}] to
conclude that the maximal ideal space of $A(K)$ coincides with $K$. (Note that
the assumption that $K$ is uniformly $H$-convex, needed in the above
approximation argument, also guarantees that $K$ is a so-called $S_\delta$ in
Rossi's terminology, i.e.\ the intersection of a sequence of pseudoconvex
domains.) If we denote the $CR$ conjugate normal by $(g_1,\ldots, g_n)$, then
by the above approximation result, $g_j\in A(K)$, $1\le j\le n$. Because the
$g_j$ have no common zeros on $K$ (see (\ref{eq3}) above), Rossi's result
implies that the ideal they generate is all of $A(K)$ (since they are not
contained in any maximal ideal). In particular, there exist $f_1,\ldots,
f_n\in A(K)$ such that $\sum\limits^n_{j=1} f_jg_j\equiv 1$ on $K$. Since
$(g_1,\ldots, g_n)$ is {\em conjugate\/} normal, this says that the (complex)
inner product of $(f_1,\ldots, f_n)$ with the normal $(\bar g_1,\ldots, \bar
g_n)$ is identically equal to 1, hence so is its real part. Consequently,
$(f_1,\ldots, f_n)$ is transverse to $b\Omega$ on $K$. Approximation of
$f_1,\ldots, f_n$ by functions holomorphic in a neighborhood of $K$ gives a
holomorphic vector field near $K$ that is transverse to $b\Omega$. This proves
b) and completes the proof of the proposition.
\end{rem}

\newpage

\centerline{REFERENCES}

\begin{enumerate}
\item\label{B-F} 
Eric Bedford and John Erik Fornaess:\ Domains with pseudoconvex neighborhood
systems, Invent.\ Math. {\bf 47} (1978), 1--27.

\item\label{B} 
Errett Bishop:\ Differentiable manifolds in complex Euclidean space, Duke
Math.\ J. {\bf 32} (1965), 1--21.

\item\label{B-S3}
Harold P.\ Boas and Emil J.\ Straube: \ Sobolev estimates for the
$\bar\partial$-Neumann operator on domains in $\mathbb{C}^n$ admitting a
defining function that is plurisubharmonic on the boundary, Math.\ Z. {\bf
206} (1991), 81--88.

\item\label{B-S1}
Harold P.\ Boas and Emil J.\ Straube: \ De Rham cohomology of manifolds
containing the points of infinite type, and Sobolev estimates for the
$\bar\partial$-Neumann problem, J.\ Geometric Analysis {\bf 3}, Nr.~3 (1993),
225--235.

\item\label{B-S2}
Harold P.\ Boas and Emil J.\ Straube: \ Global regularity of the
$\bar\partial$-Neumann problem: \ A survey of the $\mathcal{L}^2$-theory,  in
Several Complex Variables, M.\ Schneider and
Y.-T.\ Siu eds., MSRI Publications, vol.~37,  Cambridge Univ.\ Press, 1999.

\item\label{C-C} 
Alberto Candel and Lawrence Conlon:\ Foliations I, Graduate Studies in
Mathematics, vol.~23, American Mathematical Society, 2000.

\item\label{C-S}
So-Chin Chen and Mei-Chi Shaw:\ Partial Differential Equations in Several
Complex Variables, American Mathematical Society/International Press, 2001.

\item\label{DA1}
John P.\ D'Angelo:\ Iterated commutators and derivatives of the Levi form,
 in Complex Analysis, S.G.~Krantz, ed., Springer Lecture
Notes in Math., vol.~1268, Springer, 1987.

\item\label{DA2} 
John P.\ D'Angelo:\ Several Complex Variables and the Geometry of
Real Hypersurfaces, Studies in Advanced Mathematics, CRC Press, Boca Raton, FL,
1993.

\item\label{F-K}
G.B.\ Folland and J.J.\ Kohn: \ The Neumann Problem for the Cauchy-Riemann
Complex, Annals of Math.\ Studies {\bf 75}, Princeton Univ.\ Press, 1972.

\item\label{F-N}
John Erik Fornaess and Alexander Nagel: \ The Mergelyan property for weakly
pseudoconvex domains, Manuscripta Math. {\bf 22} (1977), 199--208.

\item\label{F}
Franc Forstneri\v c:\ Complex tangents of real surfaces in complex surfaces,
Duke Math.\ J. {\bf 67} (1992), 353--376.

\item\label{R}
Hugo Rossi:\ Holomorphically convex sets in several complex variables, Ann.\
Math. {\bf 24}, No.~3 (1961), 470--493.

\item\label{S}
Emil J.\ Straube: \ Good Stein neighborhood bases and regularity of the
$\bar\partial$-Neumann problem, Illinois J. Math., to appear.

\item\label{S-S}
Emil J.\ Straube and Marcel K.\ Sucheston: \ Levi foliations in pseudoconvex
boundaries and vector fields that commute approximately with $\bar\partial$,
preprint.
\end{enumerate}

\vspace{.5in}

\baselineskip = 12pt

\hspace{2.5in} Department of Mathematics

\hspace{2.5in} Texas A\&M University

\hspace{2.5in} College Station, TX \ 77843-3368

\hspace{2.5in} straube@math.tamu.edu

\end{document}